\documentclass[12pt]{amsart}

\usepackage{amsfonts,amsmath,latexsym,amssymb,verbatim,amsbsy,times}
\usepackage{amsthm}
\usepackage{pstricks}

\theoremstyle{plain}
\newtheorem{THEOREM}{Theorem}[section]

\newtheorem{theorem}[THEOREM]{Theorem}
\newtheorem{corollary}[THEOREM]{Corollary}
\newtheorem{lemma}[THEOREM]{Lemma}

\theoremstyle{definition}

\newtheorem{definition}[THEOREM]{Definition}

\theoremstyle{remark}

\newcommand{\thm}[1]{Theorem~\ref{#1}}

\newcommand{\N}{\ensuremath{\mathbb{N}}}   
\newcommand{\R}{\ensuremath{\mathbb{R}}}   

\def \RR {\mathbb{R}}
\def \a {\alpha}
\def \b {\beta}

\def \e {\epsilon}
\def \f {\varphi}

\def \l {\lambda}

\def \n {\nabla}

\def\cF{{\mathcal F}}

\def \< {\langle}
\def \> {\rangle}
\def \p {\partial}
\def \ra {\rightarrow}
\def \ss {\subset}

\DeclareMathOperator{\tr}{Tr} %

\begin{document}

\title[Solutions of 3D NSE in $B^{-1}_{\infty,\infty}$]
{On the regularity of weak solutions of the 3D Navier-Stokes
equations in $B^{-1}_{\infty,\infty}$}

\author{A. Cheskidov}
\address[A. Cheskidov]
         {Department of Mathematics\\
           University of Chicago\\
           5734 S. University Avenue\\
           Chicago, IL 60637}
\email{acheskid@math.uchicago.edu}

\author{R. Shvydkoy}
\thanks{The work of R. Shvydkoy was partially supported by NSF grant
DMS--0604050}
\address[R. Shvydkoy]
{Department of Mathematics, Stat. and Comp. Sci.\\
        University of Illinois\\
        Chicago, IL 60607}
\email{shvydkoy@math.uic.edu}

\begin{abstract}
We show that if a Leray-Hopf solution $u$ to the 3D Navier-Stokes
equation belongs to $C((0,T]; B^{-1}_{\infty,\infty})$ or its jumps
in the $B^{-1}_{\infty,\infty}$-norm do not exceed a constant
multiple of viscosity, then $u$ is regular on $(0,T]$. Our method
uses frequency local estimates on the nonlinear term, and yields an
extension of the classical Ladyzhenskaya-Prodi-Serrin criterion.
\end{abstract}


\maketitle

\section{Introduction}
We study the 3D incompressible Navier-Stokes equations (NSE)
\begin{equation} \label{NSE}
\left\{
\begin{aligned}
&\p_t u - \nu \Delta u + (u \cdot \nabla)u + \nabla p = 0, \qquad x \in \mathbb{R}^3, t\geq 0,\\
&\nabla \cdot u =0,\\
&u(0)=u_0,
\end{aligned}
\right.
\end{equation}
where $u(x,t)$, the velocity, and $p(x,t)$, the pressure, are unknowns,
$u_0 \in L^2(\mathbb{R}^3)$ is the initial condition,
and $\nu>0$ is the kinematic  viscosity coefficient of the fluid.

In 1934 Leray \cite{Leray} proved that
for every divergence-free initial data $u_0 \in L^2(\mathbb{R}^3)$
there exists a weak solution to the system \eqref{NSE} on $[0,\infty)$
with $u(0) = u_0$. Moreover, one can find a weak solution satisfying the
energy inequality (see Section~\ref{preliminaries}), called the Leray-Hopf solution.

A famous open question, equivalent to one of the millennium
problems, is whether all Leray-Hopf solutions to \eqref{NSE} on
$[0,\infty)$ with smooth initial data are regular. A weak solution
$u(t)$ is called regular if $\|u(t)\|_{H^1}$ is continuous. Even
though the regularity problem is far from been solved, numerous
regularity criteria were proved since the work of Leray. The first
result of this kind is due to Leray in the $\RR^3$ case
(see also Ladyzhenskaya, Serin, and Prodi \cite{L,S,P} in the case of a
bounded domain). It states that every
Leray-Hopf solution $u$ to \eqref{NSE} with $u \in L^{r}((0,T);
L^{s})$ is regular on $(0,T]$ provided $2/r+3/s = 1$, $s \in (3,
\infty]$. The limit case $r=\infty$ appeared to be more challenging.
First, Leray \cite{L} showed that Leray-Hopf solutions in
$L^{\infty}((0,T); L^p)$ are regular provided $p>3$. Later von Wahl
\cite{W} and Giga \cite{G} increased the regularity space to
$C((0,T]; L^{3})$, and, finally,
Esperanza, Seregin, and \v{S}ver\'ak \cite{ESS} increased it further
to $L^{\infty}((0,T); L^{3})$.

Due to this remarkable result, if a solution looses regularity, then
its $L^3$-norm blows up. The space $L^3$, as well as all the other
critical spaces for the 3D NSE (see Cannone \cite{C}), is
continuously embedded in the Besov space $B^{-1}_{\infty,\infty}$.
Hence, a natural question is whether $B^{-1}_{\infty,\infty}$-norm
has to blow up if a solution looses regularity. In this paper we
prove a slightly weaker statement. Namely, if a Leray-Hopf solution
$u$ to \eqref{NSE} is continuous in $B^{-1}_{\infty,\infty}$ on the
interval $(0,T]$, then $u$ is regular on $(0,T]$. In fact we will
prove that regularity of the solution is guaranteed even if $u$ is
discontinuous but all its jumps are smaller than a constant proportional to
the viscosity (see \thm{t:main}). Thus, discontinuities of a limited
size are allowed. As a consequence, we conclude that Leray-Hopf
solutions with $B^{-1}_{\infty,\infty}$-norm bounded by $c\nu$ are
regular. This last result will also be shown in the ''almost
everywhere" version (see Corollary \ref{c:ae}).

A related direction in the study of the 3D NSE concerns solutions
with small initial data. The best present result due to Koch and
Tataru \cite{KT} states that if the initial data is small in
$BMO^{-1}$, then there exists a global in time regular solution of
the 3D NSE with this initial data. Note that the space $BMO^{-1}$ is
continuously embedded in $B^{-1}_{\infty,\infty}$. The small
initial data theorem in this larger space is unknown, but
would be desirable to obtain.

We carry out our arguments in the case of $\R^3$, but they apply to
the periodic case as well.


\section{Preliminaries} \label{preliminaries}

\subsection{Besov spaces}

We will use the notation $\lambda_q = 2^q$ (in some inverse length
units). Let $B_r$ denote the ball centered at $0$ of radius $r$
in $\R^{3}$. Let us fix a nonnegative radial function $\chi \in {C_0^{\infty}} (B_2)$ such that $\chi(\xi)=1$ for $|\xi|\leq
1$. We further define
\begin{equation} \label{defvf}
\f(\xi) = \chi(\l_1^{-1}\xi) - \chi(\xi).
\end{equation}
and
\begin{equation} \label{defvfq}
\f_q(\xi) = \f(\l_q^{-1}\xi).
\end{equation}
For a tempered distribution vector field $u$ let us denote
\begin{align*}
u_q & = \cF^{-1}(\f_q) \ast u, \quad \text{for} \quad q > -1,\\
u_{-1} & = \cF^{-1}(\chi) \ast u,
\end{align*}
where $\cF$ denotes the Fourier transform. So, we have
$$
u = \sum_{q=-1}^\infty u_q
$$
in the sense of distributions. We also use the following notation
$$
u_{\leq q} = \sum_{p=-1}^q u_p.
$$
Let us recall the definition of Besov spaces. A tempered distribution $u$ belongs to $B^s_{p,\infty}$ for some $s\in \RR$ and $p\geq 1$ iff
$$
\|u\|_{B^s_{p,\infty}} = \sup_{q\geq -1} \l_q^s \|u_q\|_p < \infty.
$$

\subsection{Weak solutions of the 3D NSE}
In this section we recall some of the classical definitions and
results on weak solutions of the NSE.

\begin{definition}
A weak solution of \eqref{NSE} on $[0,T]$ (or $[0,\infty)$ if $T=\infty$)
is a function
$u:[0,T] \to L^2(\RR^3)$ in the class
\[
u \in  C_{\mathrm{w}}([0,T];L^2(\RR^3)) \cap L^2_{\rm loc}([0,T];H^1(\RR^3)),
\]
satisfying
\begin{multline}\label{weakform}
(u(t),\f(t)) + \int_0^t \left\{ - (u,\p_t \f) + \nu (\n u,\n \f) + (u\cdot \n u, \f) \right\} \, ds\\ = (u_0,\f(0)),
\end{multline}
for all $t\in[0,T]$ and all test functions
$\f \in C_0^\infty([0,T]\times \RR^3)$ with $\n \cdot \f = 0$. Here
$(\cdot,\cdot)$ stands for the $L^2$-inner product.
\end{definition}

\begin{theorem}[Leray] \label{thm:Leray}
For every $u_0 \in L^2(\RR^3)$, there exists a weak solution of \eqref{NSE}
on $[0,\infty)$ satisfying the following energy inequality
\begin{equation}\label{SEI}
\|u(t)\|_2^2 + 2\nu \int_{t_0}^t \|\n u(s)\|_2^2 \, ds \leq \|u(t_0)\|_2^2,
\end{equation}
for almost all $t_0 \in [0,\infty)$ including $t_0=0$, and all $t \in [t_0,\infty)$.
\end{theorem}

\begin{definition}
A Leray-Hopf solution of \eqref{NSE} on $[0,T]$ is a weak solution
on $[0,T]$ satisfying the energy inequality \eqref{SEI}
for almost all $t_0 \in (0,T)$ and all $t \in [t_0,T]$.
\end{definition}

A weak solution $u(t)$ of \eqref{NSE} on a time interval $I$ is called regular
if $\|u(t)\|_{H^1}$ is continuous on $I$.
The following is a well-known fact concerning Leray-Hopf solutions.
Note that weak solutions are not known to satisfy this property.
\begin{theorem}[Leray] \label{t:LHreg}
Let $u(t)$ be a Leray-Hopf solution of \eqref{NSE} on $[0,T]$. If
for every interval of regularity $(\alpha,\beta) \subset (0,T)$
\[
\limsup_{t \to \beta -} \|u(t)\|_{H^{s}} < \infty,
\]
for some $s > 1/2$, then $u(t)$ is regular on $(0,T]$.
\end{theorem}

\section{Main result}

In this section we state and prove our main result. Its proof is a
consequence of the limiting case of other regularity criteria in the
subcritical range of integrability parameter. We will discuss how
those criteria extend some of the classical and newly found results
on regularity.

\begin{theorem}\label{t:main}
Let $u(t)$ be a Leray-Hopf solution of \eqref{NSE} on $[0,T]$. If
$u(t)$ satisfies
\begin{equation} \label{eq:regcriterion}
\sup_{t\in(0,T]}\limsup_{t_0 \to t-} \|u(t) -
u(t_0)\|_{B^{-1}_{\infty,\infty}} < c\nu,
\end{equation}
where $c>0$ is some absolute constant, then $u(t)$ is regular on
$(0,T]$.

In particular, if $$u \in C((0,T];B^{-1}_{\infty,\infty}),$$ then
$u(t)$ is regular on $(0,T]$.
\end{theorem}

The following criterion will be a crucial ingredient in the proof of
\thm{t:main}.

\begin{lemma}\label{L:1}
Let $u(t)$ be a Leray-Hopf solution of \eqref{NSE} on $[0,T]$. There
is a constant $c>0$ such that if $u(t)$ satisfies
\begin{equation} \label{1}
\limsup_{q \to \infty} \sup_{t \in (0,T)} \l_q^{-1}
\|u_q(t)\|_{\infty} < c \nu,
\end{equation}
then $u(t)$ is regular on $(0,T]$.
\end{lemma}
\begin{proof}
Let $(\alpha, \beta) \subset (0,T]$ be an interval of regularity. Let us fix $\e\in(0,1)$ and use the weak formulation of the NSE with the test-function
$\l_q^{1+\e}(u_q)_q$. Then on $(\alpha,\beta)$ we obtain
\begin{equation}\label{nseineq-1}
\frac{1}{2}\frac{d}{dt}\l_q^{1+\e}\|u_q\|_2^2 + \nu \l_q^{3+\e}\|u_q\|_2^2
\leq \l_q^{1+\e}\int \tr[ (u\otimes u)_q \cdot \n u_q]\, dx.
\end{equation}
We can write
\[
(u \otimes u)_q = r_q(u,u) + u_q \otimes u + u \otimes u_q,
\]
for all $q>-1$, where
\[
r_q(u,u)(x) = \int_{\RR^3} \cF^{-1}(\f_q)(y) (u(x-y) - u(x)) \otimes (u(x-y) - u(x)) \, dy.
\]
After proper cancelations we arrive at
\begin{equation*}
 \int_{\RR^3} \tr[ (u\otimes u)_q \cdot \n u_q] \,dx = \int_{\RR^3} r_q(u,u) \cdot \n u_q \,dx -
 \int_{\RR^3} u_q \cdot \n u_{\leq q+1} \cdot u_q \,dx.
\end{equation*}
Let us now estimate the first term using H\"{o}lder and Bernstein's inequalities
\[
\int_{\RR^3} r_q(u,u) \cdot \n u_q \,dx \lesssim \|r_q(u,u)\|_{3/2} \l_q \|u_q\|_3.
\]
Using Littlewood-Paley decomposition we obtain as in \cite{ccfs}
\[
\begin{split}
\|r_q(u,u)\|_{3/2} &\lesssim
\int_{\RR^3} \left |\cF^{-1}(\f_q)(y) \right |\|u(\cdot - y) - u(\cdot)\|_3^2 \,dy \\
&\lesssim  \int_{\RR^3} \left |\cF^{-1}(\f_q)(y) \right| \left(
\sum_{p =-1}^q |y|^2 \l_p^2 \|u_p\|_3^2 + \sum_{p=q+1}^\infty \|u_p\|_3^2\right)  dy \\
&\lesssim \sum_{p = -1}^q  \l_q^{-2} \l_p^2 \|u_p\|_3^2 + \sum_{p=q+1}^\infty \|u_p\|_3^2.
\end{split}
\]
Analogously,
$$
\int_{\RR^3} u_q \cdot \n u_{\leq q+1} \cdot u_q \, dx \lesssim  \|u_q\|_3^2 \sum_{p=-1}^{q+1} \l_p \|u_p\|_3.
$$

Thus, up to a constant multiple independent of $u$, the nonlinear term in
\eqref{nseineq-1} is bounded by
\[
\l_q^\e \|u_q\|_3 \sum_{p =-1}^q \l_p^2 \|u_p\|_3^2 + \l_q^{2+\e}\|u_q\|_3 \sum_{p = q+1}^\infty \|u_p\|_3^2 +
\l_q^{2+\e} \|u_q\|_3^2 \sum_{p = -1}^{q+1} \l_p \|u_p\|_3.
\]
Let us now fix $Q \in \N$ and sum over $q \geq Q$ in \eqref{nseineq-1} on
$(\alpha,\beta)$ obtaining
\begin{equation} \label{eq:ineq123}
\frac{1}{2}\frac{d}{dt}\sum_{q\geq Q} \l_q^{1+\e}\|u_q\|_2^2 + \nu \sum_{q = Q}^\infty \l_q^{3+\e}\|u_q\|_2^2  \lesssim
I + II + III,
\end{equation}
where
\begin{align*}
 I & = \sum_{q = Q}^\infty \l_q^\e \|u_q\|_3 \sum_{p =-1}^q \l_p^2 \|u_p\|_3^2,\\
II & = \sum_{q = Q}^\infty\l_q^{2+\e}\|u_q\|_3 \sum_{p = q+1}^\infty \|u_p\|_3^2,\\
III& = \sum_{q = Q}^\infty \l_q^{1+\e} \|u_q\|_3^2 \sum_{p = -1}^{q+1} \l_p \|u_p\|_3.
\end{align*}
We will show that the right hand side of \eqref{eq:ineq123} obeys the following estimate
\begin{equation}\label{iii}
 I+II+III \lesssim \sum_{q = Q}^\infty \l_q^{2+\e} \|u_q\|_3^3 + R(Q),
\end{equation}
where
\[
R(Q) = \sup_{t \in (0,T)} \sum_{q =-1}^{Q-1} \l_q^{2+\e} \|u_q\|_3^3.
\]
Once this is proved, the argument goes as follows. By interpolation,
$\|u_q\|_3^3 \leq \|u_q\|_2^{2} \|u_q\|_\infty$. So,
\[
\frac{1}{2}\frac{d}{dt}\sum_{q=Q}^\infty \l_q^{1+\e}\|u_q\|_2^2 + \nu \sum_{q = Q}^\infty \l_q^{3+\e}\|u_q\|_2^2 \leq
C \sum_{q = Q}^\infty \l_q^{3+\e} \|u_q\|_2^2 ( \l_q^{-1} \|u_q\|_\infty) + CR(Q),
\]
where $C>0$ is an absolute constant  independent of $u$.
Let $c = \nu/C$. Choosing $Q$ so that $\l_q^{-1} \|u_q(t)\|_\infty < c$ for all $q \geq Q$ and all $t \in (0,T)$, we obtain
\[
\frac{1}{2}\frac{d}{dt}\sum_{q=Q}^\infty \l_q^{1+\e}\|u_q\|_2^2 \leq  CR(Q),
\]
for any $t \in (\alpha,\beta)$. Since $u(t)$ is bounded in $L^2$ on $(0,T]$, it follows that $R(Q)<\infty$. Hence the $H^{1/2+\e}$-norm of $u(t)$ is bounded on $(\a,\b)$. Therefore $u(t)$ is regular
on $(0,T]$ due to Theorem~\ref{t:LHreg}.

Let us now prove \eqref{iii}. Indeed, using Young and Jensen's inequalities
we obtain
\begin{align*}
 I & = \sum_{q = Q}^\infty \l_q^\e \|u_q\|_3 \sum_{p =-1}^q \l_p^2 \|u_p\|_3^2\\
& = \sum_{q = Q}^\infty \l_q^{\frac{2}{3} + \frac{\e}{3}} \|u_q\|_3 \sum_{p =-1}^q(\l_p\l_q^{-1})^{\frac{2}{3} - \frac{2\e}{3}} \l_p^{\frac{4}{3} + \frac{2\e}{3}} \|u_p\|_3^2 \\
& \lesssim \sum_{q = Q}^\infty \l_q^{2+\e} \|u_q\|_3^3 + \sum_{q = Q}^\infty \sum_{p =-1}^q(\l_p\l_q^{-1})^{\frac{2}{3} - \frac{2\e}{3}} \l_p^{2+\e} \|u_p\|_3^3\\
& \lesssim \sum_{q = Q}^\infty \l_q^{2+\e} \|u_q\|_3^3 + R(Q).
\end{align*}
Similarly,
\begin{align*}
 II & =  \sum_{q = Q}^\infty\l_q^{2+\e}\|u_q\|_3 \sum_{p = q+1}^\infty \|u_p\|_3^2\\
& = \sum_{q = Q}^\infty \l_q^{\frac{2}{3} + \frac{\e}{3}} \|u_q\|_3 \sum_{p = q+1}^\infty (\l_q\l_p^{-1})^{\frac{4}{3}+\frac{2\e}{3}} \l_p^{\frac{4}{3}+\frac{2\e}{3}} \|u_p\|_3^2 \\
& \lesssim \sum_{q = Q}^\infty \l_q^{2+\e} \|u_q\|_3^3 + \sum_{q = Q}^\infty \sum_{p = q+1}^\infty(\l_q\l_p^{-1})^{\frac{4}{3}+\frac{2\e}{3}} \l_p^{2+\e} \|u_p\|_3^3 \\
&  \lesssim \sum_{q = Q}^\infty \l_q^{2+\e} \|u_q\|_3^3.
\end{align*}
And finally,
\begin{align*}
 III & =  \sum_{q = Q}^\infty \l_q^{1+\e} \|u_q\|_3^2 \sum_{p = -1}^{q+1} \l_p \|u_p\|_3\\
& =  \sum_{q = Q}^\infty \l_q^{\frac{4}{3}+\frac{2\e}{3}} \|u_q\|_3^2 \sum_{p = -1}^{q+1} (\l_p\l_q^{-1})^{\frac{1}{3} - \frac{\e}{3}}  \l_p^{\frac{2}{3} + \frac{\e}{3}} \|u_p\|_3 \\
& \lesssim \sum_{q = Q}^\infty \l_q^{2+\e} \|u_q\|_3^3 + \sum_{q = Q}^\infty\sum_{p = -1}^{q+1}(\l_p\l_q^{-1})^{\frac{1}{3} - \frac{\e}{3}} \l_p^{2+\e} \|u_p\|_3^3\\
& \lesssim \sum_{q = Q}^\infty \l_q^{2+\e} \|u_q\|_3^3  + R(Q).
\end{align*}
\end{proof}

It is worth mentioning that this lemma immediately yields the
following corollary. For mild solutions in the sense of Kato a
similar result was proved in \cite{CF}.

\begin{corollary}\label{c:ae}
Let $u(t)$ be a Leray-Hopf solution of \eqref{NSE} on $[0,T]$.
There is a constant $c>0$ such that if $u(t)$ satisfies
$$\|u\|_{L^\infty((0,T);B^{-1}_{\infty,\infty})} < c \nu,$$ then
$u(t)$ is regular on $(0,T]$.
\end{corollary}

Now we are in a position to prove \thm{t:main}.

\begin{proof}[Proof of Theorem~\ref{t:main}]
Take any interval of regularity $(\alpha,\beta) \subset (0,T]$.
Since $u(t)$ is a Leray-Hopf solution, $u(t) \in
C^{\infty}(\mathbb{R}^3)$ for all $t \in (\alpha,\beta)$. Hence
\[
\phi(t):=\limsup_{q \to \infty} \l_q^{-1}\|u_q(t)\|_\infty=0, \qquad
\forall t \in (\alpha,\beta).
\]
Then, thanks to \eqref{eq:regcriterion}, $\phi(\beta) < c\nu$. This
together with \eqref{eq:regcriterion} implies that there exists $t_0
\in (\alpha, \beta)$, such that
\[
\limsup_{q \to \infty} \sup_{s \in (t_0,\beta)}
\l_q^{-1}\|u_q(s)\|_\infty < 2 c\nu.
\]
Thus, if $c$ is chosen to be half of what is in Lemma~\ref{L:1}, then that
lemma implies regularity of  $u(t)$ on $(\alpha, \beta]$, i.e.,
continuity of $\|u(t)\|_{H^1}$ on $(\alpha, \beta]$. Hence $u(t)$ is
regular on $(0,T]$ in view of Theorem~\ref{t:LHreg}.
\end{proof}

\section{Extension of the Ladyzhenskaya-Prodi-Serin criteria} \label{s:lemma:main}
We conclude with the following regularity criterion in the subcritical
range of integrability exponent.

\begin{lemma} \label{l:main}
Let $u(t)$ be a Leray-Hopf solution of \eqref{NSE} on $[0,T]$.
There is a constant $c>0$ such that if $u(t)$ satisfies
\begin{equation} \label{2}
\sup_{t \in (0,T]} \lim_{t_0 \to t-} \limsup_{q \to \infty} \int_{t_0}^t \left( \l_q^{\frac{2}{r}-1} \|u_q(s)\|_{\infty}
\right)^r \, ds< \nu^{r-1}c^r\left(\frac{r}{r-1}\right)^{r-1},
\end{equation}
for some $r\in(2,\infty)$, then $u(t)$ is regular on $(0,T]$.
\end{lemma}
\begin{proof}
Let $(\alpha, \beta) \subset (0,T]$ be an interval of regularity. We
use the weak formulation of the NSE with the test-function
$\l_q^{2r-2}(u_q)_q$. Then on $(\alpha,\beta)$ we obtain
\begin{equation}\label{nseineq}
\frac{1}{2}\frac{d}{dt}\l_q^{2r-2}\|u_q\|_2^2 \leq -\nu \l_q^{2r}\|u_q\|_2^2
+ \l_q^{2r-2}\int_{\RR^3} \tr[ (u\otimes u)_q \cdot \n u_q]\, dx.
\end{equation}
As in the proof of Lemma~\ref{L:1},
we have the following bound on the nonlinear term
\[
 \int_{\RR^3} \tr[ (u\otimes u)_q \cdot \n u_q] \,dx \leq
C_1 \sum_{p=-1}^\infty \l_{|q-p|}^{-\frac{2}{3}} \l_p \|u_p\|_3^3,
\]
where $C_1>0$ is an absolute constant.
Thanks to Young and Jensen's inequalities,
\begin{multline*}
\l_q^{2r-2}\|u_q\|^{2r-2} C_1\sum_{p=-1}^\infty \l_{|q-p|}^{-\frac{2}{3}} \l_p \|u_p\|_3^3 \leq \nu \l_q^{2r} \|u_q\|_2^{2r} \\+
(rA)^{-1} \sum_{p=-1}^\infty \l_{|q-p|}^{-\frac{2}{3}} \l_p^r \|u_p\|_3^{3r},
\end{multline*}
where
\[
A = \nu^{r-1}C^{-r}\left(\frac{r}{r-1}\right)^{r-1},
\]
and $C>0$ is an absolute constant  independent of $u$.
Therefore on $(\alpha,\beta)$ we have
\begin{equation} \label{eq:qw}
\begin{split}
\frac{1}{2r}\frac{d}{dt} \l_q^{2r-2} \|u_q\|_2^{2r} &\leq
 -\nu \l_q^{2r} \|u_q\|_2^{2r} + \l_q^{2r-2}\|u_q\|^{2r-2} C_1\sum_{p=-1}^\infty \l_{|q-p|}^{-\frac{2}{3}} \l_p \|u_p\|_3^3\\
& \leq (rA)^{-1} \sum_{p=Q+1}^\infty \l_{|q-p|}^{-\frac{2}{3}} \l_p^r \|u_p\|_3^{3r} + R(Q),
\end{split}
\end{equation}
where
\[
R(Q) = (rA)^{-1} \sup_{t \in (0,T)}
\sum_{p =-1}^{Q} \l_p^r \|u_p\|_3^{3r},
\]
$r\in(2,\infty)$, and $Q\in \mathbb{N}$. Note that since $u(t)$ is bounded in $L^2$ on $(0,T]$,
we have $R(Q)<\infty$.

Now assume that \eqref{2} holds for $c=(C\sqrt{24})^{-1}$. Then there exist
$t_0 \in (\alpha, \beta)$ and $Q$ large enough, such that
\[
\sup_{q>Q} \int_{t_0}^\beta(\l_q^{\frac{2}{r}-1}\|u_q(s)\|_\infty)^r \, ds \leq
{\textstyle \frac{1}{24}} A.
\]
Let $t\in(t_0,\beta)$. Integrating \eqref{eq:qw} we obtain
\begin{multline} \label{eq:ghjk}
\sup_{s\in(t_0,t)}\l_q^{2r-2} \|u_q(s)\|_2^{2r} -
\l_q^{2r-2} \|u_q(t_0)\|_2^{2r}\\ \leq
2A^{-1} \int_{t_0}^t \sum_{p=Q+1}^\infty  \l_{|q-p|}^{-\frac{2}{3}} \l_p^r \|u_p(s)\|_3^{3r} \, ds +R(Q)(t-t_0).
\end{multline}
Thanks to Levi's convergence theorem
we have
\[
\begin{split}
\int_{t_0}^t \sum_{p=Q+1}^\infty & \l_{|q-p|}^{-\frac{2}{3}} \l_p^r \|u_p(s)\|_3^{3r} \, ds\\
&\leq \sup_{p> Q} \left\{ \sup_{s\in(t_0,t)}\l_p^{2r-2} \|u_p(s)\|_2^{2r}\int_{t_0}^t(\l_p^{\frac{2}{r}-1}\|u_p(s)\|_\infty)^r \, ds
 \right\} \sum_{p=1}^\infty  \l_{|q-p|}^{-\frac{2}{3}}\\
&\leq 6\sup_{p> Q} \left\{\sup_{s\in(t_0,t)}\l_p^{2r-2} \|u_p(s)\|_2^{2r} \right\}\sup_{p> Q} \int_{t_0}^t(\l_q^{\frac{2}{r}-1}\|u_p(s)\|_\infty)^r \, ds
\\
&\leq {\textstyle \frac{1}{4}}A  \sup_{p> Q} \sup_{s\in(t_0,t)}\l_p^{2r-2} \|u_p(s)\|_2^{2r}.
\end{split}
\]
Hence, taking the supremum of \eqref{eq:ghjk} over $q>Q$, we obtain
\[
{\textstyle \frac{1}{2}}\sup_{q> Q} \sup_{s\in(t_0,t)}\l_q^{2r-2} \|u_q(s)\|_2^{2r} \leq
\sup_{q> Q} \l_q^{2r-2} \|u_q(t_0)\|_2^{2r} + R(Q)(t-t_0).
\]
Letting $t \ra \b$ we now have
$$
\sup_{q > Q} \sup_{s \in (t_0, \b)} \l_q^{1- \frac{1}{r}} \|u_q(s)\|_2 < \infty.
$$
Thus, since $r>2$, the solution $u(t)$ is bounded in $H^{s}$
on $(t_0, \beta)$ for some $s>1/2$, which implies regularity
in view of \thm{t:LHreg}.
\end{proof}

Let us now see how this lemma
implies an extension of the classical Ladyzhenskaya-Serin-Prodi
condition in the subcritical range of integrability exponent.
Indeed, by Bernstein inequalities and Littlewood-Paley Theorem we
have the inclusions
\begin{align*}
B^{\frac{3}{s} + \frac{2}{r}-1}_{s,\infty} & \ss
B^{\frac{2}{r}-1}_{\infty,\infty},\\
L^s &\ss B^0_{s,\infty},
\end{align*}
for all $s\geq 1$.
On the other hand, as an easy consequence of Lemma~\ref{l:main} we
have the following theorem.
\begin{theorem}\label{c:crit}
Let $u(t)$ be a Leray-Hopf solution of \eqref{NSE} on $[0,T]$
satisfying
\begin{equation}\label{e:crit}
u \in L^{r}((0,T); B^{\frac{2}{r}-1}_{\infty,\infty}),
\end{equation}
for some $r \in (2,\infty)$. Then $u(t)$ is regular on $(0,T]$.
\end{theorem}
Let us also note that for a negative smoothness parameter as in
\eqref{e:crit}, the homogeneous version of the Besov space is smaller
than the nonhomogeneous one, i.e.,
$$
\dot{B}^{\frac{2}{r}-1}_{\infty,\infty} \ss
B^{\frac{2}{r}-1}_{\infty,\infty}, \quad \text{for all} \quad r >2.
$$
In view of this fact, Theorem~\ref{c:crit} extends the
corresponding result obtained recently by Q. Chen and Z. Zhang in
\cite{CZ}.

\bibliographystyle{plain}

\end{document}